\documentclass[10 pt, leqno]{amsart}

\usepackage {amsfonts}
\usepackage{amsthm}
\usepackage{amssymb}
\usepackage{latexsym}
\usepackage{amsmath}
\usepackage{mathrsfs}

\theoremstyle{plain}
\newtheorem{tw}{Theorem}[section]

\newtheorem {lem} [tw]{Lemma}
\newtheorem {prop}[tw] {Proposition}

\newtheorem{cor}[tw]{Corollary}

\newtheorem {exam}[tw] {Example}

\theoremstyle{definition}
\newtheorem {deft}[tw] {Definition}

\newcommand{\bc} {\mathbb C}
\newcommand{\Oo} {\mathcal O}
\newcommand{\mJ} {\mathcal J}
\newcommand{\bn}{\mathbb N}

\newcommand{\be}{\mathbb E}

\newcommand{\bz}{\mathbb Z}

\newcommand{\comp}{\mathcal{K}(E)}
\newcommand{\adj}{\mathcal{L}(E)}

\newcommand{\alg} {\mathsf{A}}

\newcommand {\Ran} {{\textrm{Ran}}}
\newcommand {\id} {{\textrm{id}}}

\newcommand{\Hil}{\mathsf{H}}

\newcommand{\Kil}{\mathsf{K}}

\newcommand{\Nk}{\bn_0^k}
\newcommand{\Njk}{\bn_0^{k \setminus j}}
\newcommand{\Nalk}{\bn_0^{k, \alpha}}
\newcommand{\Nbek}{\bn_0^{k, \beta}}
\newcommand{\wT}{\wt{T}}
\newcommand{\llin}{\mathcal{L}}
\newcommand{\FFock}{\mathcal{F}}
\newcommand{\onek}{\{1,\ldots,k\}}
\newcommand{\wS}{\wt{S}}
\newcommand{\wV}{\wt{V}}

\newenvironment{rlist}
{

\begin{enumerate}}
{\end{enumerate}}

\newcommand{\la}{\langle}
\newcommand{\ra}{\rangle}
\newcommand{\indfe}{\pi^{\FFock(\be)}}

\newcommand{\ot}{\otimes}

\newcommand{\wt}{\widetilde}

\numberwithin{equation}{section}


\begin{document}
\author{Adam Skalski}
\author{Joachim Zacharias}
\footnote{\emph{Permanent address of the first named author:}
Department of Mathematics, University of \L\'{o}d\'{z}, ul.
Banacha 22, 90-238 \L\'{o}d\'{z}, Poland.

Research supported by the EPSRC grant no. RIS 24893}

\address{School of Mathematical Sciences,  University of Nottingham,
Nottingham, NG7 2RD} \email{adam.skalski@nottingham.ac.uk} \email{joachim.zacharias@nottingham.ac.uk  }

\title{\bf Wold decomposition for representations of product systems of $C^*$-correspondences}
\begin{abstract}
\noindent Higher-rank versions of Wold decomposition are shown to hold for doubly commuting isometric
representations of product systems of $C^*$-correspondences over $\Nk$, generalising the classical result for a doubly commuting pair of isometries due to M.\,S\l oci\'nski.
Certain decompositions are also obtained for the general, not necessarily doubly commuting, case and several
corollaries and examples are provided. Possibilities of extending isometric representations to fully coisometric ones are discussed.
\end{abstract}

\maketitle

The classical notion of Wold decomposition refers to the unique decomposition of a Hilbert space isometry into
a part which is unitary and a part which is isomorphic to a unilateral shift. For a simple proof and several
applications of this result we refer to the classical monograph \cite{Nagy}. In analogy with the famous dilation
problem for tuples of contractions, it is natural to ask whether some version of Wold decomposition is available
for a tuple of commuting isometries. Indeed, M.\,S\l oci\'nski established in \cite{Sloc} such a decomposition for
a \emph{doubly commuting} pair of isometries. This result (and its generalisations) was later used in
\cite{Coburn} to provide models for tuples of commuting isometries and analyse the structure of $C^*$-algebras they
generate. Another example of the analysis of the structure of a pair of commuting isometries, also of relevance to
our work here, may be found in \cite{Popo}.

In recent years there has been an increased interest in Wold decompositions for objects of a different type. It
originated from the work of G.\,Popescu, who in \cite{Popescu} established a result of this kind for a row
contraction. Various related ideas were extended to an impressive degree in the series of papers of P.\,Muhly and
B.\,Solel, who developed the theory of tensor algebras over $C^*$-correspondences. In particular in \cite{MSWold}
they proved the existence of a Wold decomposition for an isometric representation of a $C^*$-correspondence over
a $C^*$-algebra $\alg$. Another, more concrete, example of such a decomposition may be found in \cite{jury}.

In this paper we establish a higher-rank version of the main result of \cite{MSWold}. The corresponding crucial concept of a product system of a $C^*$-correspondence over $\Nk$ was introduced in \cite{fowl}; the notion has been exploited in the recent work by B.\,Solel on dilations of commuting completely positive maps ([So$_{1-2}$]). Here we prove that every doubly commuting isometric representation of a product system of $C^*$-correspondences over $\Nk$ decomposes uniquely into a combination of \emph{fully coisometric} and \emph{induced} parts (in the classical terminology they correspond respectively to a unitary and a shift part). It turns out that for isometric representations which are not doubly commuting it is still possible to characterise maximal pieces of the carrier
space which have analogous properties; we also obtain an appropriate generalisation of the results of \cite{Popo}.
The main motivation behind this work is the hope that it will be of use for understanding the structure of
$C^*$-algebras generated by commuting tuples of row isometries, thus taking further the ideas of \cite{Coburn}.
The class above is easily seen to contain tensor products of Cuntz algebras.

The plan of the paper is as follows: the first section introduces the notation, recalls the basic results of
P.\,Muhly and B.\,Solel and contains a short proof of the existence of maximal fully coisometric summands for
isometric representations (of a product system). Section 2 is devoted to doubly commuting representations, with
the Wold decomposition established and characterisation of those contractive representations whose minimal
isometric dilations are induced or fully coisometric. In Section 3 the assumption of double commutativity is
dropped, and the analogue of the decomposition in \cite{Popo} obtained. Section 4 provides examples of the
structure provided by our Wold decompositions for some concrete classes of product systems of
$C^*$-correspondences. Finally in Section 5 we discuss the possibility of constructing fully coisometric extensions
of isometric representations, with certain positive results obtained for the case of a single $C^*$-correspondence.

\section{Basic notation 
and the existence of  maximal fully coisometric summands}

Let $\alg$ be a $C^*$-algebra. By a $C^*$-correspondence $E$ over $\alg$ is meant a Hilbert $C^*$-module over
$\alg$, equipped with the structure of a left $\alg$-module (via a nonzero $*$-homomorphism $\phi$ mapping
$\alg$ into the $C^*$-algebra of adjointable operators on $E$). We always assume that
$E$ is \emph{essential} as a left $\alg$-module, i.e.\ the closed linear span of $\phi(\alg) E$ is equal to $E$. Each $C^*$-correspondence is considered with the
usual operator space structure (i.e.\ the one coming from viewing them as corners in respective linking algebras). In the
following we will frequently use the internal tensor product construction for $C^*$-Hilbert modules -- for the
details of this and other aspects of the theory of the latter we refer to \cite{Lance}.
All representations of $\alg$ in this paper are nondegenerate. 

\begin{deft}
Let $\Hil$ be a Hilbert space. By a (completely contractive covariant) representation of a $C^*$-correspondence
$E$ on $\Hil$ is meant a pair $(\sigma, T)$, where $(\sigma,\Hil)$ is a  representation of $\alg$ on $\Hil$,
and $T:E \to B(\Hil)$ is a linear completely contractive map such that
\[ T(a \xi b) = \sigma(a) T(\xi) \sigma(b), \;\;\; a,b \in \alg, \xi \in E.\]
It is called isometric if for each $\xi, \eta \in E$
\[ T(\xi)^* T(\eta) = \sigma(\la \xi, \eta \ra).\]
\end{deft}

A representation $(T, \sigma)$ determines a contraction $\wT : E \ot_{\sigma} \Hil \to \Hil$ given by $\wT (\xi
\ot h) = T(\xi) h$ ($\xi \in E, h \in \Hil$). This satisfies:
\begin{equation} \label{tilde} \wT (\phi(a) \ot I_{\Hil}) = \sigma(a) \wT, \;\;\; a \in \alg
\end{equation}
($\phi$ denoting the left action of $\alg$ on $E$), and one can in fact show that, given a representation $\sigma$,
there is a 1-1 correspondence
between contractions satisfying \eqref{tilde} and representations of $E$ (\cite{MSgen} Lemma 2.1). The isometric
representations are exactly those for which $\wT$ is an isometry. A (completely contractive covariant)
representation $(T, \sigma)$ is called \emph{fully coisometric} if $\wT \wT^* = I_{\Hil}$.

Fix for the rest of the paper $k \in \bn$ and write $\bn_0 =\bn\cup\{0\}$. The basic object for this
work will be a
\emph{product system $\be$ of $C^*$-correspondences over $\bn_0^k$} (\cite{fowl}). As explained in [So$_{1-2}$],
$\be$ can be thought of as a family of $k$ $C^*$-correspondences $\{E_1, \ldots,E_k\}$ together with the unitary
isomorphisms $t_{i,j}: E_i \ot E_j \to E_j \ot E_i$ ($i>j$). This point of view entails identifying for all
$n=(n_1, \ldots, n_k) \in \Nk$ the correspondence $\be (n)$ with $E_1^{\ot^{ n_1}} \ot \cdots \ot E_k^{\ot^{
n_k}}$. We additionally write $t_{i,i} = \id_{E_i \ot E_i}$, $t_{i,j} = t_{j,i}^{-1}$ for $i<j$, and denote the
`basis' elements of $\Nk$ by $e_1,\ldots,e_k$.

\begin{deft}
Let $\be$ be a product system over $\Nk$. By a (covariant completely contractive) representation of $\be$ on a
Hilbert space $\Hil$ is meant a tuple $(\sigma, T^{(1)}, \ldots, T^{(k)})$, where $(\sigma,\Hil)$ is a
 representation of $\alg$, $T^{(i)}:E_i \to B(\Hil)$ are linear completely contractive maps such that
\[ T^{(i)}(a \xi_i b) = \sigma(a) T^{(i)}(\xi_i) \sigma(b), \;\;\; a,b \in \alg, \xi_i \in E_i,\]
and
\begin{equation} \label{rep} \wT^{(i)} (I_{E_i} \ot \wT^{(j)}) = \wT^{(j)} (I_{E_j} \ot \wT^{(i)}) (t_{i,j} \ot I_{\Hil})\end{equation}
 for $i,j\in \onek$.
 Such a representation is called isometric if each $(\sigma, T^{(i)})$ is isometric as a representation
 of $E_i$, and fully coisometric if each $(\sigma, T^{(i)})$ is fully coisometric.
\end{deft}

Two representations $(\sigma, T^{(1)}, \ldots, T^{(k)})$, $(\rho, S^{(1)}, \ldots, S^{(k)})$ of $\be$,
respectively on Hilbert spaces $\Hil$ and $\Kil$, are called isomorphic if there exists a unitary $U:\Hil \to
\Kil$ implementing the unitary equivalence of representations $\sigma$ and $\rho$ and such that for all $i \in
\onek, \xi \in E_i$ there is $S^{(i)}(\xi) = U T^{(i)} (\xi) U^*$.

Let $(\sigma, T^{(1)}, \ldots, T^{(k)})$ be an isometric representation of $\be$ on $\Hil$. For each $i \in \onek$
and $l \in \bn$ define $\wT^{(i)}_l: E_i^{\ot l} \ot_{\sigma} \Hil\to \Hil$ by the formula
\[ \wT^{(i)}_l (\xi_1 \ot \cdots \ot \xi_l \ot h) = T^{(i)} (\xi_1) \cdots T^{(i)}(\xi_l) h\]
($\xi_1, \ldots, \xi_l \in E_{i}, h \in \Hil$). It is easy to see that each $\wT^{(i)}_l$ is isometric, with the
range equal to the closed linear span of $\{T^{(i)} (\xi_1) \cdots T^{(i)}(\xi_l)h: \xi_1, \ldots \xi_l \in E_{i},
h \in \Hil\}$. We will write $P_l^i$ for the orthogonal projection on this set, so that $P^i_l =  \wT^{(i)}_l
\wT^{(i)^*}_l$. The following is Lemma 2.3 of \cite{MSWold} (in our notation):

\begin{lem} \label{endomor}
Let $i \in \onek$. The formula
\[L_i(x) = \wT^{(i)} (I_{E_i} \ot x) (\wT^{(i)^*}), \;\;\; x \in \sigma(\alg)'\]
defines a normal endomorphism of the commutant of $\sigma(\alg)$. Moreover
\[L_i^l(x) = \wT^{(i)}_l (I_{E_i} \ot x) \wT^{(i)^*}_l, \;\;\; x \in \sigma(\alg)',\]
$L_i^l(I) = P_l^i$ ($l \in \bn$).
\end{lem}

The lemma above uses only the fact that each $(\sigma, T^{(i)})$ is an isometric representation of $E^i$. In our
context a simple calculation shows that
\begin{align*} L_i(L_j(x)) =&\,  \wT^{(i)} (I_{E_i} \ot (\wT^{(j)} (I_{E_j} \ot x) \wT^{(j)^*}))
\wT^{(i)^*} \\
=&\,  \wT^{(i)} (I_{E_i} \ot \wT^{(j)}) ( I_{E_i \ot E_j} \ot x) (\wT^{(i)} (I_{E_i} \ot \wT^{(j)}))^*
\\
=&\, \wT^{(j)} (I_{E_j} \ot \wT^{(i)})  (t_{i,j} \ot I_{\Hil}) ( I_{E_i \ot E_j} \ot x)
   (t_{i,j}^* \ot I_{\Hil}) (\wT^{(j)} (I_{E_j} \ot \wT^{(i)}))^* \\=&\,
   L_j(L_i(x)).\end{align*}

In particular, we obtain the action of $\Nk$ on $\sigma(\alg)'$ by normal endomorphisms:
 $L(n) := L_1^{n_1} \cdots L_k^{n_k}$, $n \in \Nk$ (with $L(0) =I_{\sigma(\alg)'}$).

\begin{deft}
A subspace $\Kil \subset \Hil$ is called reducing for the representation $(\sigma, T^{(1)}, \ldots, T^{(k)})$ if
it reduces $\sigma(\alg)$ (so that the projection onto $\Kil$, denoted further by $P_{\Kil}$, lies  in
$\sigma(\alg)')$, and both $\Kil, \Kil^{\perp}$ are left invariant by all operators $T^{(i)} (\xi_i)$ for $i \in
\onek$, $\xi_i \in E_i$. Then it is easy to see that the obvious `restriction' procedure yields a representation
of $\be$ on $\Kil$, which is called a summand of $(\sigma, T^{(1)}, \ldots, T^{(k)})$ and will be
denoted by $(\sigma, T^{(1)}, \ldots, T^{(k)})|_{\Kil}$.
\end{deft}

The following is a straightforward corollary of Lemma 2.5 of \cite{MSWold}:

\begin{cor}  \label{reduce}
Let $\Kil$ be a Hilbert subspace of $\Hil$. It reduces $(\sigma, T^{(1)}, \ldots, T^{(k)})$ if and only if it
reduces $\sigma(\alg)$ and $L_i(P_{\Kil}) = P_{\Kil}P^i_1 =  P^i_1 P_{\Kil}$  for each $i \in \onek$.
\end{cor}

Let $P^i_{\infty}$ denote the infimum (equivalently, the limit) of the sequence of projections $\{P_l^i:l \in
\bn\}$. The arguments above show that
\[ \lim_{n \in \Nk} L(n) (I) = \lim_{l \to \infty} (L_1\cdots L_k)^l (I) =
\inf_{n \in \Nk} L(n) (I);\] the projection given by the formula above will be denoted $P_{\infty}$. Note that we
obviously have $P_{\infty} \leq  \bigwedge_{i=1}^k P_{\infty}^i$; is there equality here in general? We suspect
not, it should be possible to construct a counterexample already via models of pairs of commuting contractions
considered in \cite{Coburn}.

Now we are ready to establish the existence of maximal fully coisometric summands for any representation of $\be$.
The proposition below will in fact turn out to be a special case of Theorem \ref{main2}; we give a separate proof
here as we believe that this case helps to establish the intuition for what follows.

\begin{prop} \label{unitary}
Every isometric representation $(\sigma, T^{(1)}, \ldots, T^{(k)})$ of $\be$ on a Hilbert space $\Hil$ has a
unique maximal fully coisometric summand.
\end{prop}

\begin{proof}
We will show that the required summand is given by $\Hil_{\infty}:=P_{\infty}\Hil$. It is obvious that
$L_i(P_{\infty}) = P_{\infty}$, and as $P_\infty \leq P^i_{1}$, we have $L_i(P_{\infty}) =
P_{\infty}P^i_{1}=P^i_{1}P_{\infty}$. The fact that the summand given by restriction to $\Hil_{\infty}$ is fully
coisometric follows easily: if we denote the restricted representation by $(\sigma', S^{(1)}, \ldots, S^{(k)})$,
then for each $i \in \onek$ $\wt{S}^{(i)} = P_{\infty} \wT^{(i)} (I_{E_i} \ot P_{\infty})$, and
\[ \wt{S}^{(i)} \wt{S}^{(i)^*} = P_{\infty} \wT^{(i)} (I_{E_i} \ot P_{\infty}) \wT^{(i)^*}
 P_{\infty}=
P_{\infty} L_i ( P_{\infty}) P_{\infty} = P_{\infty}.\] It remains to show maximality (uniqueness will follow).
For this assume that $P$ is a projection on a reducing subspace of our representation, such that the resulting
restriction is fully coisometric. The first condition yields $L_i (P) = P_1^i P = P P_1^i $, the second $P L_i(P)
P = P$ (again, $i \in \onek$). This implies that $L_i(P) = P$, so $P =(L(n)) (P)\leq
L(n) (I)$ for arbitrary $n \in \Nk$, and $P \leq P_{\infty}$.
\end{proof}

Note that $\Hil_{\infty}$ is in fact equal to the maximal fully coisometric summand appearing in the Wold
decomposition of the isometric representation $T^{(1)} \ot \cdots \ot T^{(k)}$ of the $C^*$-correspondence $E_1
\ot \cdots \ot E_k$.

\section{Doubly commuting case}

Recall the definition of the doubly commuting representation of $\be$ introduced in \cite{Solk}:

\begin{deft}  \label{dcom}
A representation $(\sigma, T^{(1)}, \ldots, T^{(k)})$ of  $\be$ on a Hilbert space $\Hil$ is called doubly
commuting if for each $i,j \in \{1,\ldots,k\}$, $i \neq j$ implies
\begin{equation}\wT^{(j)^*} \wT^{(i)} =
 (I_{E_j} \ot \wT^{(i)})  (t_{i,j} \ot I_{\Hil})  (I_{E_i} \ot \wT^{(j)^*}).
 \label{doubly}\end{equation}
\end{deft}

 B.\,Solel showed in \cite{Solk} that for isometric representations the condition \eqref{doubly} is equivalent to \emph{Nica-covariance} (\cite{Nica}). It is immediate that summands of doubly commuting representations are doubly commuting. Note also that
if $(\sigma, T^{(1)}, \ldots, T^{(k)})$ is a representation of $\be$ on $\Hil$ and for some $i\in \onek$ $(\sigma,
T^{(i)})$ is fully coisometric, then the condition \eqref{doubly} holds for arbitrary $j \in \onek
\setminus\{i\}$.

Note the crucial property following from the formula \eqref{doubly}: this time the family of projections $\{P_l^i:
i \in \onek, l \in \bn\}$ defined before Lemma \ref{endomor} is commutative. We will in fact show more; if we
define
\begin{equation} P(0) = I_{\Hil}, \;\;\;  P(n)=(L(n))(I) \label{Pn}\end{equation}
($n \in \Nk\setminus\{0\}$), then
\begin{equation}  \label{crucial} P(m) P(n) = P(m \vee n)\end{equation}
($m,n \in \Nk$). The last formula can be proved inductively, starting first with  ($i, j \in \onek$, $i \neq j$)
\begin{align*}P_1^i  P_1^j =& \,\wT^{(i)}  \wT^{(i)^*}   \wT^{(j)}  \wT^{(j)^*} =
\wT^{(i)}  (I_{E_i} \ot \wT^{(j)})  (t_{j,i} \ot I_{\Hil})  (I_{E_j} \ot \wT^{(i)^*})
(\wT^{(j)})^* \\
=& \,\wT^{(i)} (I_{E_i} \ot \wT^{(j)}) (t_{j,i} \ot I_{\Hil})
  (\wT^{(i)} (I_{E_i} \ot \wT^{(j)}) (t_{j,i} \ot I_{\Hil}))^* \\
  =& \,\wT^{(i)} (I_{E_i} \ot \wT^{(j)}) (t_{j,i} \ot I_{\Hil}) (t_{i,j} \ot I_{\Hil})
(I_{E_i} \ot \wT^{(j)^*}) \wT^{(i)^*} = L_i (P_1^j).\end{align*} In the next step we show that for $i, j$ as
above, $l \in \bn$
\begin{align*} P_1^i  P_{l+1}^j =&\, \wT^{(i)} \wT^{(i)^*}  \wT^{(j)} (I_{E_j} \ot P_{l}^j) \wT^{(j)^*} =
   \wT^{(j)} (I_{E_j} \ot P_1^i)  (I_{E_j} \ot P_{l}^j) \wT^{(j)^*} \\
=&\,\wT^{(j)}   ( I_{E_j} \ot P(e_i + l e_j))  \wT^{(j)^*} =
   L_j (P(e_i + l e_j)) = P(e_i+ (l+1) e_j).\end{align*}
Then, in similar vein, for $i, j, l $ as above and $p \in \bn$
\[ P_p^i  P_l^j =  P(p e_i +l e_j).\]
Further for $ i \in \onek$, $l, p \in \bn$, $ l \leq p$
\[ P_p^i P_l^i = L_i^p(I_{\Hil}) L_i^l (I_{\Hil}) = L_i^l (L_i^{p-l}(I_{\Hil})) L_i^l (I_{\Hil})
=  L_i^l (L_i^{p-l}(I_{\Hil})) =  L_i^p(I_{\Hil}) = P_{l\vee p}^i.\]
Combining the formulas above, and once again reasoning inductively 
yields finally the equality \eqref{crucial}.
 This implies in particular that  $\{P(n): n \in \Nk\}$ is a lattice. It can be easily extended to a complete lattice
by allowing also the relevant `limit' projections, using the indexing set $(\bn_0 \cup \{\infty\})^k$ instead of
$\Nk$.

\vspace*{0.2cm}

\noindent \textbf{Induced representations (Fock space shifts).}

\vspace*{0.1 cm}

\noindent Assume that $F$ is a $C^*$-correspondence over $\alg$, and that $\pi$ is a representation of $\alg$ on a
Hilbert space $\Kil$. This induces a representation $\pi^{F}$ of $\llin (F)$ on the Hilbert space $F \ot_{\pi}
\Kil$, given by the formula $\pi^{F}(T) = T \ot I_{\Kil}$ ($T \in \llin(F)$).

Let now $\be$ be a product system of $C^*$-correspondences over $\bn_0^k$ and let $\FFock(\be)$ denote the Fock
module of $\be$ (see \cite{fowl} for the details). Then the procedure given above yields the induced
representation $\pi^{\FFock(\be)}: \llin(\FFock(\be)) \to \FFock(\be) \ot_{\pi} \Kil$. This in turn determines a
representation of $\be$ by the following formulas:
\[ \sigma(a) = \indfe (\phi_{\infty} (a)) = \phi_{\infty} (a) \ot I_{\Kil}, \;\;\; a \in \alg,\]
\[ T^{(i)} (\xi_i) = \indfe (T_{\xi_i}) = T_{\xi_i} \ot I_{\Kil}, \;\;\; i \in \onek, \xi_i \in E_i,\]
where $\phi_{\infty}$ denotes the canonical left action of $\alg$ on $\FFock(\be)$ and $T_{\xi_i}$ denotes a
Toeplitz creation operator determined by $\xi_i$. Any representation of $\be$ constructed in the way described
above is called an \emph{induced representation}.  It may be thought of as a generalised shift; it reduces to the
usual $k$-shift with multiplicity $\textrm{dim } \Kil$ when $\alg = \bc$ and $\be$ is given by the trivial family
$(\bc, \ldots, \bc)$. Note that it is nondegenerate due to the assumption that all $C^*$-correspondences are essential
as left $\alg$-modules.

The so-obtained representation of $\be$ may be easily shown to be isometric. It is also doubly
commuting; this can be observed if we note that for induced representations the maps $\wt{T}^{(i)}$ reduce to
formal identifications coming from coassociativity of various tensor products. For a formal proof we refer to
\cite{fowl}.

With the formula \eqref{crucial} in hand, we can easily prove the following higher-rank version of Corollary 2.10
of \cite{MSWold}:

\begin{lem} \label{induced}
Let $\be$ be a product system of $C^*$-correspondences over $\bn_0^k$ and let \\$(\sigma, T^{(1)}, \ldots,
T^{(k)})$ be a doubly commuting isometric representation of $\be$ on a Hilbert space $\Hil$. Then $(\sigma,
T^{(1)}, \ldots, T^{(k)})$ is isomorphic to an induced representation if and only if $P^i_{\infty} = 0$ for each
$i \in \onek$, where $P^i_{\infty}$ are projections introduced before Proposition \ref{unitary}.
\end{lem}

\begin{proof}
The `only if' direction follows easily by inspecting the Fock space situation (think about the analogy with the
shift model).

For `if' we will again use the notation established in \eqref{Pn}. Denote for $n \in \Nk$
\[Q(n) = (P_{n_1}^1 - P_{n_1+1}^1) \cdots (P_{n_k}^k - P_{n_k+1}^k).\]
Note that due to \eqref{crucial} one can equivalently write
\[Q(n) = (P(n) - P(n + e_1)) \cdots (P(n) - P(n + e_k)).\]
It is easy to see, using the first formula, that the projections from the family $\{Q(n): n \in \Nk\}$ are
mutually orthogonal; moreover the assumption on vanishing of the limits $P^i_{\infty}$ yields $\bigoplus_{n \in
\Nk} Q(n) \Hil = \Hil$. Let $\Hil_0:= Q(0) \Hil$ (this will be our `wandering subspace') and let $\sigma_0=
\sigma|_{\Hil_0}$. Note that $L(n) (Q(0)) = Q(n)$ (it follows from the second formula describing $Q(n)$). This
in particular implies that if we write
\begin{equation}
\wT(n) = \wT^{(1)}_{n_1} \cdots (I_{\be(n_1 e_1 + \cdots + n_{k-1} e_{k-1})} \ot \wT^{(k)}_{n_k}): \be(n) \ot \Hil
\to \Hil
 \label{Tn} \end{equation}
then $\wT(n) (I_{\be(n)} \ot Q(0))$ is a partial isometry ($E\ot_{\sigma} \Hil \to \Hil$)  with the initial
projection $I_{\be(n)} \ot Q(0)$  and the range projection $Q(n)$.

Recall that $\FFock(\be) = \bigoplus_{n \in \Nk} \be(n)$ and define the map $U: \FFock(\be)\ot_{\sigma_0} \Hil_0
\to \Hil$ by
\[ U \left((h_n)_{n \in \Nk}\right) =  \sum_{n \in \Nk} \wT(n) h_n\]
($(h_n)_{n \in \Nk} \in \bigoplus_{n \in \Nk} \be(n) \ot_{\sigma_0} \Hil_0$). It remains to see that due to the
remarks above $U$ is a Hilbert space isomorphism, and moreover it implements the isomorphism of $(\sigma, T^{(1)},
\ldots, T^{(k)})$ with the representation induced by $\sigma_0$. The last statement may be checked directly by
applying the definition of operators $\wT^{(i)}$ and exploiting the formula \eqref{tilde} in conjunction with
properties of internal tensor products.
\end{proof}

Identically as in Proposition 2.11 of \cite{MSWold}, we can deduce the following

\begin{cor} \label{subinduced}
If $(\sigma, T^{(1)}, \ldots, T^{(k)})$ is a representation of $\be$ on $\Hil$, isomorphic to an induced
representation, and $\Kil\subset \Hil$ is a reducing subspace, then 
$(\sigma, T^{(1)}, \ldots, T^{(k)})|_{\Kil}$ is also isomorphic to an induced representation.
\end{cor}

\vspace*{0.2cm}

\noindent \textbf{Decomposition for doubly commuting representations.}

\vspace*{0.1 cm}

The following theorem is the main result of this section. It generalises (and in a sense puts in the same context)
the Wold decompositions from \cite{Sloc} and  \cite{MSWold}.

\begin{tw} \label{main}
Every doubly commuting isometric representation $(\sigma, T^{(1)}, \ldots, T^{(k)})$ of $\be$ on a Hilbert space
$\Hil$ has a unique decomposition given by $\Hil = \bigoplus_{\alpha \subset \{1,\ldots,k\}} \Hil_{\alpha}$ such
that for each $\alpha \subset \{1,\ldots,k\}$, $\alpha=\{\alpha_1, \ldots, \alpha_r\}$,
\begin{rlist}
\item $\Hil_{\alpha}$ reduces
$(\sigma, T^{(1)}, \ldots, T^{(k)})$;
\item $(\sigma, T^{(\alpha_1)}, \ldots,  T^{(\alpha_r)})|_{\Hil_{\alpha}}$ is isomorphic to
an induced representation of the product subsystem $\be_{\alpha}$ over $\bn_0^r$, given by the
$C^*$-correspondences $E_{\alpha_1}, \ldots, E_{\alpha_r}$ (and obvious isomorphisms $t_{\alpha_i, \alpha_j})$;
\item for each $i \in \{1,\ldots,k\} \setminus \alpha$ the representation
$(\sigma, T^{(i)})|_{\Hil_{\alpha}}$ of $E_i$ is fully coisometric.
\end{rlist}
\end{tw}

\begin{proof}
Note that $\alpha=\emptyset$ corresponds to the summand described in Proposition \ref{unitary}. The decomposition
will be given by $\Hil_{\alpha} = P_{\alpha} \Hil$, where for each $\alpha  \subset \{1,\ldots,k\}$,
$\alpha=\{\alpha_1, \ldots, \alpha_r\}$, $\{1,\ldots,k\} \setminus \alpha:=\{\beta_1, \ldots, \beta_{k-r}\}$,
\[P_{\alpha} = P_{\infty}^{{\alpha_1}^{\perp}} \cdots P_{\infty}^{{\alpha_r}^{\perp}}
P_{\infty}^{\beta_1} \cdots P_{\infty}^{\beta_{k-r}}.\] Fix for now $\alpha$ as above. As ($i,j \in \onek$)
\[L_i(P^j_{\infty}) = \lim_{l \in \bn} P(e_i + l e_j) = \lim_{l \in \bn}  P^j_l P^i_1 =
  P^j_{\infty} P^i_1,\]
and similarly $L_i(P^j_{\infty}) = P^i_1 P^j_{\infty} $, we automatically obtain
\[ L_i (P_{\alpha}) = P_{\alpha} P^i_1 =  P^i_1  P_{\alpha},\]
so from Corollary \ref{reduce} it follows that $\Hil_{\alpha}$ is reducing for our representation.
To prove the condition (ii) it is enough (by Lemma \ref{induced}) to show that for each $i \in \{1,\ldots,r\}$
\[ P_{\alpha} P_{\infty}^{\alpha_i} P_{\alpha} = 0.\]
This follows immediately from the definition of $P_{\alpha}$. Further if $j \in \{1, \ldots, k-r\}$, we obviously
have
\[P_{\alpha} L_{\beta_j} (P_{\alpha}) P_{\alpha} =
 P_{\alpha} P_1^{\beta_j} P_{\alpha} = P_{\alpha},\]
 so the condition (iii) is satisfied.

It remains to show the uniqueness of the decomposition satisfying the conditions in the theorem. Suppose that
$\Hil = \bigoplus_{\alpha \subset \{1,\ldots,k\}} \Hil'_{\alpha}$ for some family $\{\Hil'_{\alpha}:\alpha \subset
\{1,\ldots,k\}\}$ of subspaces of $\Hil$ satisfying  (i)-(iii) above. Fix for a moment $\alpha$ and let
$P_{\alpha}'$ denote the projection on  $\Hil'_{\alpha}$. As $\Hil'_{\alpha}$ is reducing, there is ($i \in
\onek$)
\[ L_i(P_{\alpha}') = P_{\alpha}' P_1^i = P_1^i P_{\alpha}' =
  P_{\alpha}' P_1^i  P_{\alpha}'.\]
As for $i\in \onek \setminus \alpha$ the representation $(\sigma, T^{(i)})|_{\Hil_{\alpha}'}$ is fully
coisometric, there is
\[ P_{\alpha}' L_i (P_{\alpha}')  P_{\alpha}' = P_{\alpha}'.\]
As in Proposition \ref{unitary}, we deduce from the above that $P_{\alpha}' \leq P_{\infty}^i$. Finally by Lemma
\ref{induced} if $j \in \alpha$, then
\[P_{\alpha}'  P_{\infty}^j  P_{\alpha}' = 0,\]
so that $P_{\alpha}' \leq I_{\Hil} - P_{\infty}^j$. This in turn yields $P_{\alpha}' \leq P_{\alpha}$, and as
$\bigoplus_{\alpha \subset \{1,\ldots,k\}} \Hil_{\alpha} = \bigoplus_{\alpha \subset \{1,\ldots,k\}}
\Hil'_{\alpha}$, we must have $\Hil_{\alpha} = \Hil_{\alpha}'$ for all $\alpha \subset \onek$.
\end{proof}

Observe that in the doubly commuting case the equality $\inf_{n \in \Nk} P(n) = \bigwedge_{i=1}^k P_{\infty}^i$
holds automatically.

\vspace*{0.2cm}

\noindent \textbf{Remarks on connections with minimal isometric dilations of contractive representations.}

\vspace*{0.1 cm} In Theorem 3.5 of \cite{Solk} B.\,Solel proved that every doubly commuting completely
contractive covariant representation
 $(\sigma, T^{(1)}, \ldots, T^{(k)})$ of $\be$ has a (unique up to a unitary equivalence) \emph{minimal regular dilation}
 to an isometric
representation. This means that there exists an isometric representation $(\sigma', V^{(1)}, \ldots, V^{(k)})$ of
 $\be$ on a Hilbert space $\Kil$ such that $\Hil \subset \Kil$, each $V^{(i)} (\xi_i)$ ($i \in \onek, \xi_i \in E_i$)
 leaves $\Hil$ invariant, $\Hil$ is reducing for $\sigma'$ with $\sigma'|_{\Hil} = \sigma$, for all $m,n \in \Nk$
\[ (I_{\be(m)} \ot P_{\Hil}) \wV(m)^* \wV(n)|_{\be(n) \ot_{\sigma} \Hil} = (I_{\be(m)} \ot P_{\Hil}) \wT(m)^* \wT(n)\]
and $\Kil$ contains no nontrivial subspaces containing $\Hil$ and invariant for each $V^{(i)} (\xi_i)$ ($i \in
\onek, \xi_i \in E_i$). The dilated representation may be also shown to be doubly commuting. Using this fact we
can obtain the following generalisations of Propositions 2.3 and 2.5 of \cite{Popescu}.

\begin{prop} \label{dilatafter}
Let $(\sigma, T^{(1)}, \ldots, T^{(k)})$  be a doubly commuting completely contractive covariant representation of
$\be$ on a Hilbert space $\Hil$. Its minimal isometric dilation $(\sigma', V^{(1)}, \ldots, V^{(k)})$ is
isomorphic to an induced representation if and only if for all $h \in \Hil$, $j\in \onek$
\begin{equation} \label{nounit} \lim_{l \to \infty} \|(\wT(le_j))^* h\| = 0.\end{equation}
\end{prop}
\begin{proof}
If $(\sigma', V^{(1)}, \ldots, V^{(k)})$ is an induced representation, then its carrier Hilbert space is
isomorphic to  $\FFock(\be) \ot_{\pi} \Kil'$ for some representation $(\pi, \Kil')$ of $\alg$. Moreover for all $k=(k_n)_{n \in \Nk} \in \FFock(\be) \ot_{\pi} \Kil'$, $j \in
\onek$, $l \in \bn$
\[\wV(le_j)^* k = (k_n)_{n \geq le_j} \in \be(le_j) \ot \FFock(\be) \ot_{\pi} \Kil',\]
and as for $h \in \Hil$, $m \in \Nk$ there is $\wV(m)^* h = \wT(m)^* h$, condition \eqref{nounit} holds.

Assume now that condition \eqref{nounit} is satisfied. We claim that for all $m \in \Nk$, $\xi \in \be(m)$, $h'
\in \Hil$ the following holds: \begin{equation} \label{limit} \lim_{l \to \infty} \|(\wV(le_j))^* \wV(m) (\xi \ot
h')\| =0.
\end{equation}
 Indeed, for $j \in \Nk$, $l \geq m_j$
there is (by double commutativity) \begin{align*} \wV(l e_j)^* \wV(m) (\xi \ot h') &= (I_{\be(l e_j)} \ot \wV
(m-m_je_j) ) (t\ot I_{\Hil}) (\xi \ot \wV((l-m_j)e_j)^* h')\\
&=(I_{\be(l e_j)} \ot \wV (m-m_je_j) ) (t\ot I_{\Hil}) (\xi \ot \wT((l-m_j)e_j)^* h') \end{align*}
 where $t$ is the isomorphism between $\be(m) \ot
\be((l-m_j) e_j)$ and $\be(le_j) \ot \be(m-m_j e_j)$. The last formula implies \eqref{limit}.

By minimality, the linear combinations of vectors of the type $\wV(m) (\xi \ot h')$ are dense in $\Kil$. It
follows that  for all $w \in \Kil$, $j \in \onek$
\[ \lim_{l \to \infty} \|(\wV(le_j))^* w\|=0\]
As a result, all the projections $P_{\infty}^j$ from the proof of Theorem \ref{main2} are trivial. This completes the argument.
\end{proof}

\begin{prop}
Let $(\sigma, T^{(1)}, \ldots, T^{(k)})$  be a doubly commuting completely contractive covariant representation of
$\be$ on a Hilbert space $\Hil$. Its minimal isometric dilation $(\sigma', V^{(1)}, \ldots, V^{(k)})$ is fully
coisometric if and only if $(\sigma, T^{(1)}, \ldots, T^{(k)})$ is fully coisometric.
\end{prop}

\begin{proof}
The only if implication follows immediately from the fact that if $\wV^{(j)} (\wV^{(j)})^* h = h,$ for each $h \in
\Hil$ then $ \wV^{(j)} (\wV^{(j)})^* h = P_{\Hil} \wV^{(j)} (\wV^{(j)})^* h$ and the latter is  equal to $
\wT^{(j)} (\wT^{(j)})^* h$.

Assume then that $(\sigma, T^{(1)}, \ldots, T^{(k)})$ is fully coisometric and suppose that
\\$(\sigma', V^{(1)}, \ldots, V^{(k)})$ is not. There exists then $j \in \onek$ such that
$\wV^{(j)} (\wV^{(j)})^* \Kil \neq \Kil$ (where $\Kil$ denotes the carrier space of the dilated representation).
By minimality there must exist $m\in \Nk$, $\xi \in \be(m)$, $h \in \Hil$ such that
\[ \wV^{(j)} (\wV^{(j)})^* \wV(m) (\xi \ot h) \neq  \wV(m) (\xi \ot h).\]
Note that if $m_j \neq 0$, the above inequality cannot hold, as then $\wV^{(j)} (\wV^{(j)})^* \wV(m) = \wV(m)$. On the
other hand if $m_j =0$, then using double commutativity we see that
\begin{align*}\wV^{(j)} (\wV^{(j)})^* \wV(m) (\xi \ot h) &=
\wV(m) (I_{\be(m)} \ot \wV^{(j)} (\wV^{(j)})^*) (\xi \ot h) \\
&= \wV(m) (\xi \ot \wV^{(j)} (\wV^{(j)})^* h).\end{align*}
 The inequality above implies now that $\wV^{(j)}
(\wV^{(j)})^* h \neq h$; this in turn forces the existence of a nonzero $h'\in \Hil$ such that $\wV^{(j)}
(\wV^{(j)})^* h'=0$. But then $\wT^{(j)} (\wT^{(j)})^* h'=0$ and we obtain the contradiction.
\end{proof}

\section{Beyond the doubly commuting case}

Consider again an isometric, but not necessarily doubly commuting, representation of $\be$ on a Hilbert space
$\Hil$.  Although in general one cannot expect the existence of a neat decomposition analogous to the one
formulated in Theorem \ref{main}, we can still identify for each $\alpha \subset \{1,\ldots,k\}$,
$\alpha=\{\alpha_1, \ldots, \alpha_r\}$ a maximal summand of the representation such that the corresponding
restriction satisfies conditions (ii) and (iii) from Theorem \ref{main}. It will have become clear that the main
theorem below in a sense contains Theorem \ref{main}. The choice of this order of presentation is motivated by the
fact that especially the proof of maximality uses the properties of the decomposition for doubly commuting
representations established in Theorem \ref{main}.

Recall the definition \eqref{Tn} and introduce the following notation: $\Njk= \{m \in \Nk: m_j=0\}$ ($j\in
\{1,\ldots,k\}$). Moreover for $\alpha \subset\onek$ write $\Nalk=\{m\in \Nk:m_i=0 \textrm{ for } i \notin
\alpha\}$.

\begin{tw} \label{main2}
Let $(\sigma, T^{(1)}, \ldots, T^{(k)})$ be an isometric representation of $\be$ on a Hilbert space $\Hil$ and let
$\alpha \subset \{1,\ldots,k\}$, $\alpha=\{\alpha_1, \ldots, \alpha_r\}$. There exists a (unique) maximal subspace
$\Hil_{\alpha}$ of $\Hil$ among all subspaces $\Kil$ verifying the following:
\begin{rlist}
\item $\Kil$ reduces
$(\sigma, T^{(1)}, \ldots, T^{(k)})$;
\item $(\sigma, T^{(\alpha_1)}, \ldots,  T^{(\alpha_r)})|_{\Kil}$ is isomorphic to
an induced representation of the product subsystem $\be_{\alpha}$ over $\bn_0^r$, given by the
$C^*$-correspondences $E_{\alpha_1}, \ldots, E_{\alpha_r}$ (and obvious isomorphisms $t_{\alpha_i, \alpha_j})$;
\item for each $i \in \{1,\ldots,k\} \setminus \alpha$ the representation
$(\sigma, T^{(i)})|_{\Kil}$ of $E_i$ is fully coisometric.
\end{rlist}
\end{tw}

\begin{proof}
Write $\beta=\onek \setminus \alpha= \{\beta_1,\ldots,\beta_{k-r}\}$. Let for each $j\in \alpha$
\begin{equation} \Kil^{(j)} = \{h \in \Hil: \forall_{m \in \Njk} \; \forall_{\xi \in \be(m)} \;\; \xi \ot h = \wt{T}(m)^* (P_1^{j})^{\perp} \wT(m) (\xi \ot h)   \}. \label{Kj}\end{equation}
Note that
\[ \Kil^{(j)} = \{h \in \Hil: \forall_{m \in \Njk} \; \forall_{\xi \in \be(m)} \;\;  \wt{T}(m)^* P_1^{j} \wT(m)
(\xi \ot h) = 0\}   \]\[=\{h \in \Hil: \forall_{m \in \Njk} \; \forall_{\xi \in \be(m)} \;\ P_1^{j} \wT(m) (\xi
\ot h) = 0\}.\]
The latter formulas will prove more useful. Denote the orthogonal projection on
$\Kil^{(j)}$ by $R^{(j)}$. Observe that for any $m \in \Njk$ we have
\begin{equation} L(m) (R^{(j)}) \leq (P_1^{j})^{\perp}.\label{orthog}\end{equation}
To see that this holds write the projection on the left as $\wT(m) (I_{\be(m)} \ot R^{(j)})\wT(m)^*$ and note that
$P_1^{j} \wT(m) (I_{\be(m)} \ot R^{(j)}) = 0$. To compare the construction here with the one in Theorem
\ref{main} note that if the representation $(\sigma, T^{(1)}, \ldots, T^{(k)})$ is doubly commuting, $R^{(j)}$
coincides with $(P_1^{j})^{\perp}$.

If $a \in \alg$ and $h \in \Kil^{(j)}$, then $\sigma(a) h \in \Kil^{(j)}$ (this statement basically reduces to the
observation that $P_1^{j} \in \sigma(\alg)'$). This in particular implies  that $R^{(j)}$ is in $\sigma(\alg)'$.
Further denote the infimum of all $R^{(j)}$ by $R$, so that $R$ is a projection and
\[ R = \bigwedge_{j\in \alpha} R^{(j)} .\]
(By this we mean the projection onto the intersection of the $\Kil^{(j)}$.) Put $R_{\infty}=\lim_{n \in \Nbek} L(n) (R)$ and let $D(m) = L (m) (R_{\infty})$ for each $m\in \Nalk$.

We will start by showing that the projections $D(m)$ are mutually orthogonal; this is essentially the statement
that $R_{\infty}$ is a `wandering subspace' for $(T^{(\alpha_1)}, \ldots, T^{(\alpha_r)})$. Choose $m, m' \in
\Nalk$ so that $m\neq m'$ and let $j\in \alpha$ be such  that $m_j \neq m'_j$. Then we have $ D(m) =  L_{j}^{m_j}
(L(m-m_j e_j)(R_{\infty}))$ and for all $p\in \Nbek$
\[L(m-m_je_j) (L(p) (R)) = L(m-m_je_j + p) (R) \leq  L(m-m_je_j + p) (R^{(j)}) \leq (P_1^{j})^{\perp},\]
where the last inequality follows from \eqref{orthog}. This means in particular that $L(m-m_je_j) (R_{\infty})
\leq (P_1^{j})^{\perp}$, which yields
\[ D(m) \leq L_j^{m_j} ( (P_1^{j})^{\perp});\]
similarly we obtain
\[ D(m') \leq L_j^{m'_j} ((P_1^{j})^{\perp}).\]
The projections on the right hand side in two expressions above are mutually orthogonal (see the argument in Lemma
\ref{induced}, applied for $k=1$).

Let  $D= \sum_{m \in \Nalk } D(m)$ and put  $\Hil_{\alpha} = \Ran D = \bigoplus_{m \in \Nalk} D(m) \Hil$. We will
show that
 $\Hil_{\alpha}$ reduces $(\sigma, T^{(1)}, \ldots, T^{(k)})$.
To this end observe first that, for $i \in \beta$ and $m\in \Nalk$, $L_i(D(m)) = L_i(L(m) (R_{\infty})) =
L(m)(L_i(R_{\infty})) = L(m) (R_{\infty}) = D(m)$, so $L_i (D) =D$ (note that this immediately shows that the
condition (iii) holds). If $j \in \alpha$, and $m\in \Nalk \setminus \Njk$, then $D(m)=L_j (L(m-e_j)(R_{\infty}))
\leq L_j(I)$, and $P_1^{j} D(m) P_1^{j} = D(m)$. Assume then that $ m\in \Nalk \cap \Njk$. Choose $p\in \Nbek$
and note that we have
\[L(m) (L(p) (R)) = L(m+p)(R) \leq L(m+p)(R^{(j)}) \leq (P_1^{j})^{\perp}.\]
Taking the ultraweak limit (in $p$) yields
\[D(m) = L(m) (R_{\infty}) \leq (P_1^{j})^{\perp},\]
and finally we have
\begin{align*} L_j(D) =L_j \left(\sum_{m \in \Nalk} D(m)\right) & \\ = \sum_{m \in \Nalk \setminus \Njk} D(m) &=
   P_1^{j} \left(\sum_{m \in \Nalk} D(m) \right)   P_1^{j} =  P_1^{j} D    P_1^{j}.\end{align*}

The fact that $(\sigma, T^{(\alpha_1)}, \ldots, T^{(\alpha_r)})|_{\Hil_{\alpha}}$ is isomorphic to an induced representation is
immediate (an argument analogous to the one in the proof of Lemma \ref{induced} suffices).

It remains to show that $\Hil_{\alpha}$ is a maximal subspace of $\Hil$ satisfying conditions (i)-(iii). To this
end suppose that $\Kil$ is a subspace of $\Hil$ such that (i)-(iii) hold; denote the projection on $\Kil$ by $Q$.
Further denote the restricted representation $(\sigma, T^{(1)}, \ldots, T^{(k)})|_{\Kil}$ by $(\sigma', S^{(1)},
\ldots, S^{(k)})$, writing also  $S(n)= QT(n)  (I_{\be(n)} \ot Q)$ for $n \in \Nk$. Note that in particular
$(\sigma', S^{(1)}, \ldots, S^{(k)})$ is doubly commuting - this follows from the remarks after Definition
\ref{dcom} and before Lemma \ref{induced}. Then there is a decomposition
\[\Kil  = \bigoplus_{n \in \Nalk} L(n) (Q(0))\Hil,\]
where $Q(0)= Q (P_1^{\alpha_1})^{\perp}\cdots (P_1^{\alpha_r})^{\perp}Q $. To conclude the proof it is enough
to show  that $Q(0) \leq R_{\infty}$.

To this end note first that $Q(0)\leq R^{(j)}$ for each $j\in \alpha$. Indeed, if $h \in Q(0)\Hil$, $m \in \Njk$ and
$\xi \in \be(m)$ then, as $Q$ reduces the representation we started with,
\begin{align*}  \wT(m)^* P_1^{j} \wT(m) (\xi \ot h) & = \wT(m)^* P_1^{j} \wT(m) (I_{\be(m)} \ot Q(0))(\xi \ot h)\\
&=\wS(m)^* \wS(e_j) \wS(e_j)^* \wS(m) (I_{\be(m)} \ot Q(0)) (\xi \ot h).\end{align*} As the restricted
representation is doubly commuting,
\begin{align*} \wT(m)^* P_1^{j} & \wT(m)  (\xi \ot h) \\
& = (I_{\be(m)} \ot \wS(e_j)) \wS(m)^*  \wS(m) (I_{\be(m)} \ot \wS(e_j)^*) (I_{\be(m)} \ot Q(0)) (\xi \ot h) \\
& =(I_{\be(m)} \ot \wS(e_j) \wS(e_j)^*) (I_{\be(m)} \ot Q(0)) (\xi \ot h)=0.\end{align*} This implies that $Q(0)
\leq R$, so also for all $p\in \Nbek$
\[ I_{\be(p)} \ot Q(0) \leq I_{\be(p)} \ot R.\]
The projection on the left is equal to the projection $\wS(p)^* Q(0) \wS(p)$ - this follows again from the
double-commutation relation and the definition of $Q(0)$. We obtain then
\[  \wS(p)^* Q(0) \wS(p) \leq I_{\be(p)} \ot R,\]
so also
\[  \wT(p) \wS(p)^* Q(0) \wS(p) \wT(p)^*  \leq \wT(p) (I_{\be(p)} \ot R) \wT(p)^*.\]
But, as each $(\sigma,S_i)$ for  $i \in \beta$ is fully coisometric, and $p \in \Nbek$, the relation above can be
rewritten as
\[ Q(0) \leq L(p) (R).\]
 The fact that $p\in \Nbek$ was arbitrary yields the desired inequality $Q(0) \leq R_{\infty}$ and ends the proof.
\end{proof}

Note that $\alpha=\emptyset$ corresponds to the purely unitary part constructed in Lemma \ref{unitary}.

It is easy to see that if the representation $(\sigma, T^{(1)}, \ldots, T^{(k)})$ is doubly commuting, the spaces
$\Hil_{\alpha}$ in Theorem \ref{main2} coincide with the ones constructed in Theorem \ref{main}. Simple argument
yields therefore the following:

\begin{cor} \label{maxdc}
Let $(\sigma, T^{(1)}, \ldots, T^{(k)})$ be an isometric representation of $\be$ on a Hilbert space $\Hil$. Let
for each $\alpha \subset \onek$ $\Hil_{\alpha}$ denote the subspace constructed in Theorem \ref{main2}. Then the
space $\Hil_{dc} = \bigoplus_{\alpha \subset \onek} \Hil_{\alpha}$ is the maximal reducing subspace for $(\sigma,
T^{(1)}, \ldots, T^{(k)})$ such that the corresponding summand
 of $(\sigma, T^{(1)}, \ldots, T^{(k)})$ is doubly commuting.
\end{cor}

\vspace*{0.2cm}

\noindent \textbf{Weakly induced representations}

\vspace*{0.1 cm} When $k=2$ and $\be_1=\be_2=\bc$, $t_{2,1}(\mu) = \mu \; (\mu \in \bc)$, the isometric
representation of the resulting product system over $\bn_0^2$ is naturally given by a commuting pair of
isometries. The Wold decompositions for such a pair were investigated recently in \cite{Popo}. In that paper
D.\,Popovici introduced the notion of a weak bi-shift and showed that in general one can decompose the Hilbert
space into the four following components: `purely unitary' (corresponding in the language of Theorem \ref{main2} to
$\alpha=\emptyset$), `shift-unitary' ($\alpha=\{1\}$), `unitary-shift' ($\alpha=\{2\}$) and \emph{weak bi-shift}.
Analogous decompositions can be achieved in our context for isometric representations of an arbitrary
product system
of $C^*$-correspondences over $\bn_0^2$. The precise formulation of this statement and the sketch of possible
generalisations for $k >2$ are presented below.

Let $(\sigma, T^{(1)}, T^{(2)})$ be an isometric representation of $\be$  (a product system over $\bn_0^2$) on a
Hilbert space $\Hil$. Recall the definition \eqref{Kj}. Note that although in general $K^{(1)}$ does not reduce
$T^{(2)}$, it is easy to check that all operators $T^{(2)} (\xi)$, $\xi \in E_2$ leave $K^{(1)}$ invariant, so
that we can consider the isometric representation of $E_2$ on $K^{(1)}$ obtained by the restriction of $(\sigma,
T^{(2)})$ (recall that $K^{(1)}$ reduces $\sigma$). This paves the way for the following definition:

\begin{deft} \label{weak}
A representation $(\sigma, T^{(1)}, T^{(2)})$ is weakly bi-induced if each of  $(\sigma,T^{(1)} \ot T^{(2)})$,
$(\sigma,T^{(1)})|_{\Kil^{(2)}}$, $(\sigma,T^{(2)})|_{\Kil^{(1)}}$ is isomorphic to an induced representation
(respectively of $E_1 \ot E_2$, $E_1$, $E_2$).
\end{deft}

Note that if a representation is weakly bi-induced, then by analysing the explicit descriptions from the proof
Theorem \ref{main2} one can show that the maximal parts corresponding to fully coisometric, fully
coisometric--induced and induced--fully coisometric summands are trivial.

\begin{tw}
Let $(\sigma, T^{(1)}, T^{(2)})$ be an isometric representation of $\be$  (a product system over $\bn_0^2$) on a
Hilbert space $\Hil$. Then $\Hil$ has a unique decomposition into four subspaces:
\[ \Hil = \Hil_{\emptyset} \oplus \Hil_1 \oplus \Hil_2 \oplus \Hil_{wi},\]
such that
\begin{rlist}
\item each of the four subspaces is reducing for $(\sigma, T^{(1)}, T^{(2)})$;
\item $(\sigma, T^{(1)})$ is fully coisometric when restricted to $\Hil_{\emptyset}$ or $\Hil_2$, $(\sigma, T^{(2)})$ is fully coisometric when restricted to $\Hil_{\emptyset}$ or $\Hil_1$;
\item  $(\sigma, T^{(1)})|_{\Hil_1}$, $(\sigma, T^{(2)})|_{\Hil_2}$ are isomorphic to  induced representations;
\item the representation $(\sigma, T^{(1)}, T^{(2)})|_{\Hil_{wi}}$ is weakly-bi induced.
\end{rlist}
\end{tw}
\begin{proof}
Denote by $\Hil_{\emptyset}$, $\Hil_1$ and $\Hil_2$ the spaces corresponding to the subsets $\emptyset, \{1\},
\{2\} \subset \{1,2\}$ and described in Theorem \ref{main2}. Let $\Hil' = \Hil \ominus (\Hil_{\emptyset} \oplus
\Hil_1 \oplus \Hil_2)$. Observe first that $\Hil'$ reduces $(\sigma, T^{(1)}, T^{(2)})$. We will show that
$(\sigma, T^{(1)}, T^{(2)})|_{\Hil'}$ is weakly bi-induced.

Suppose this is not the case, and that, for example, $(\sigma, T^{(1)})|_{\Kil'^{(2)}}$ is not isomorphic to an induced
representation, where
\[ \Kil'^{(2)} = \{h \in \Hil': \forall_{l\in \bn_0} \; \forall_{\xi \in E_1^{\ot l}} \;\;
 P_1^{2} \wT^{(1)}_l (\xi \ot h) = 0\}.\]
Let then $Q$ be a nonzero projection corresponding to the maximal fully coisometric summand of $(\sigma,
T^{(1)})|_{\Kil'^{(2)}}$. As $\Kil'^{(2)} = \Kil^{(2)} \cap \Hil'$, $Q\leq R$, where $R$ is the projection from
the proof of Theorem \ref{main2} (the case $\alpha=\{2\}$). As the restricted representation is fully coisometric,
$L_1(Q)=Q$. This means that in particular for all $l \in \bn_0$
\[Q = L_1^l (Q) \leq L_1^l (R),\]
so also $Q\leq R_{\infty}\leq D$ (again we borrow the notation from the proof of Theorem \ref{main2},
$\alpha=\{2\}$). But this means that $Q\Hil \subset \Hil_2$, so $Q$ has to be $0$, as $Q\Hil \subset \Hil' \subset
\Hil_2^{\perp}$. The other cases may be reduced to contradiction in a similar way.

It remains to show the uniqueness. To this end assume that there is another decomposition $\Hil =
\Hil'_{\emptyset} \oplus \Hil'_1 \oplus \Hil'_2 \oplus \Hil'_{wi}$ satisfying the requirements of the theorem. By
the maximality in Theorem \ref{main2} we obtain immediately that $\Hil'_{\emptyset} \subset \Hil_{\emptyset}$,
$\Hil'_1 \subset \Hil_1$ and $\Hil'_2 \subset \Hil_2$. Now
\[ \Hil'_{wi} \supset (\Hil_{\emptyset} \ominus \Hil'_{\emptyset}) \oplus
(\Hil_1 \ominus \Hil'_1) \oplus (\Hil_2 \ominus \Hil'_2).\] By Corollary \ref{subinduced} $\Hil_1 \ominus \Hil'_1$
is an induced (with respect to the first variable)--fully coisometric (with respect to the second variable)
summand; by the remark after Definition \ref{weak} it has to be equal to $\{0\}$. Similarly $\Hil_2=\Hil'_2$,
$\Hil_{\emptyset}=\Hil_{\emptyset}$, and the proof is finished.
\end{proof}

The decomposition in the above theorem can be extended for $k>2$; for example if $k=3$ one obtains
\[\Hil = \Hil_{\emptyset}\oplus \Hil_1 \oplus \Hil_2 \oplus \Hil_3 \oplus \Hil_{1,2} \oplus
\Hil_{1,3} \oplus \Hil_{2,3} \oplus \Hil_{wi},\] where `numbered' summands correspond to the ones constructed in
Theorem \ref{main2}, and the restriction to $\Hil_{wi}$ is a weakly tri-induced representation --- that means that
all $(\sigma,T^{(1)} \ot T^{(2)} \ot T^{(3)})$, $(\sigma,T^{(1)} \ot T^{(2)})|_{\Kil^{(3)}}$, $(\sigma, T^{(1)}
\ot T^{(3)})|_{\Kil^{(2)}}$, $(\sigma,T^{(2)} \ot T^{(3)})|_{\Kil^{(1)}}$, $(\sigma, T^{(1)})|_{\Kil^{(2)}\wedge
\Kil^{(3)}}$, $(\sigma,T^{(2)})|_{\Kil^{(1)}\wedge \Kil^{(3)}}$, $(\sigma,T^{(3)})|_{\Kil^{(1)}\wedge \Kil^{(2)}}$
are isomorphic to induced representations.

These results may be extended further for arbitrary $k$, with the inductive proofs, not requiring introducing any
new methods or concepts.

\section{Examples} \label{Examples}
In this section we follow closely Section 4 of \cite{Soltwo}, where the contractive covariant representations were
given interpretations for various product systems $\be$. Here we explain what being fully coisometric means and
how induced representations look like in each of these cases.

\subsection*{Case $k=1$}

Each (completely contractive covariant) representation of a product system over $\bn_0$ corresponds, by the very
definition that we use here, to a (completely covariant contractive) representation of a single
$C^*$-correspondence; each such representation is obviously doubly commuting (this being an empty condition) and
the results of Section 2 reduce to the case thoroughly analysed in \cite{MSWold}.

\subsection*{Case $\alg = \bc$}

Here the product system of $C^*$-correspondences reduces to a product system of Hilbert spaces; the unitary
isomorphisms $t_{i,j}$ are Hilbert space unitaries from $\Hil_i \ot \Hil_j$ to $\Hil_j \ot \Hil_i$ (where we write
$\Hil_i$ instead of $E_i$). The essentiality assumptions imply that the left actions of $\bc$ on each $\Hil_i$
are simply given by multiplications.
Representations of $\be$ are given by a carrier space $\Hil$ and a set of contractions
$(\wT^{(1)}, \ldots, \wT^{(k)})$, where each $\wT^{(i)}:\Hil_i \ot \Hil \to \Hil$ and the condition \eqref{rep} is
satisfied. Note that we may think of $\wT^{(i)}$ as a row contraction: when  an orthonormal basis
$\{e_{\alpha}^{(i)}: \alpha \in \mathcal{J}_i\}$ is chosen  in each $\Hil_i$, we are precisely in a higher-rank
version of the situation studied in the series of papers of G.\,Popescu. Defining $S_{\alpha}^i (h) = \wT^{(i)}
(e_{\alpha}^{(i)} \ot h)$, $h \in \Hil, \alpha\in \mathcal{J}_i$ we obtain for each $i \in \onek$ a family of
operators on $\Hil$ (indexed by $\alpha$) such that
\begin{equation}\sum_{\alpha \in \mathcal{J}_i} S_{\alpha}^i (S_\alpha^i)^* \leq I_{\Hil}
\label{strange0}\end{equation} and the following condition is satisfied: for all $i, j \in \onek$, $\alpha \in
\mJ_i, \beta \in \mJ_j$
\begin{equation} S^i_{\alpha} S^j_{\beta} = \sum_{\alpha' \in \mJ_i, \beta' \in \mJ_j}
\la e_{\beta'}  \ot e_{\alpha'}, t_{i,j} (e_{\alpha} \ot e_{\beta}) \ra  S^j_{\beta'} S^i_{\alpha'}
\label{strange}\end{equation} (some indices $i,j$ have been dropped). The convergence above is understood
strongly.

Conversely, each family of operators $\{S_{\alpha}^i: i \in \onek, \alpha \in \mJ_i\}$ satisfying the conditions
\eqref{strange0} and \eqref{strange} yields a representation of $\be$ (it is enough to define $\wT^{(i)}$ by the
formula $\wT^{(i)} (\xi \ot h) = \sum_{\alpha \in \mJ_i} \la e_{\alpha}, \xi\ra S_{\alpha}^i (h)$, $i \in \onek$,
$\xi \in E_i$, $h \in \Hil$).

The fact that the representation $(\wT^{(1)}, \ldots, \wT^{(k)})$ is isometric corresponds to the fact that for
each $i \in \onek$ the family $\{S_{\alpha}^i: i \in \onek, \alpha \in \mJ_i\}$ consists of isometries with
orthogonal ranges; it is fully coisometric if the ranges of these isometries sum to identity on $\Hil$. Further
$(\wT^{(1)}, \ldots, \wT^{(k)})$ is doubly commuting if and only if for all $i, j \in \onek$, $\alpha \in \mJ_i,
\beta \in \mJ_j$
\begin{equation} (S^j_{\beta})^* S^i_{\alpha} = \sum_{\alpha' \in \mJ_i, \beta' \in \mJ_j}
\la e_{\beta}  \ot e_{\alpha'}, t_{i,j} (e_{\alpha} \ot e_{\beta'}) \ra  S^i_{\alpha'} (S^j_{\beta'})^*.
\label{strange1}\end{equation}

Unsurprisingly, induced representations are particularly easy to describe if each $\Hil_i$ is one-dimensional.
Note first that then $t_{i,j}$ are complex numbers of modulus 1; general isometric representation is given by $k$
isometries $S_1,\ldots, S_k$ on a Hilbert space $\Hil$ such that for all $i,j \in \onek$
\begin{equation} S_i S_j = t_{i,j} S_j S_i.\label{strange2}\end{equation}
Corollary \ref{maxdc} assures the existence of a maximal subspace of $\Hil$ reducing for all the isometries in
question and such that the reduced isometries satisfy the condition
\begin{equation} S_j^* S_i = t_{i,j} S_i S_j^*.\label{strange3}\end{equation}

The induced representation is in this case a twist of the usual shift. It is unique up to  multiplicity: if the
latter is equal to 1, $\Hil = l^2(\Nk)$
and the isometries $S_1,\ldots, S_k$ are given by the formula
\[ S_i (\delta_m) = t_{i,1}^{m_1} t_{i,2}^{m_2} \cdots t_{i,i-1}^{m_{i-1}} \delta_{m+e_i},\]
$i \in \onek$, $m \in \Nk$. It is elementary to check that $S_1,\ldots, S_k$ satisfy both \eqref{strange2} and
\eqref{strange3}. The general induced representation is obtained by tensoring $ l^2(\Nk)$ with an auxilliary
Hilbert space and ampliating isometries defined above by the identity.

Theorem \ref{main2} implies that whenever we are given isometries $S_1,\ldots, S_k$ satisfying the conditions
\eqref{strange2} and \eqref{strange3} on a Hilbert space $\Hil$, the space $\Hil$ decomposes uniquely into
subspaces reducing some of the isometries to unitaries and the other to a tuple unitarily equivalent to the twisted
higher-dimensional shift described above.

\subsection*{Case $E_i=_{\alpha_i}\alg$}
Assume we are given a $C^*$-algebra $\alg$ and $\alpha_1,\ldots, \alpha_k \in \textrm{Aut}(\alg)$. The
$C^*$-correspondences  $E_i:=_{\alpha_i}\alg$ are defined in the usual way, with the right action given by the
right multiplication, usual $\alg$-valued scalar product $ \la a,b\ra = a^*b$ ($a,b \in \alg$) and the left action
of $a\in \alg$ given by a left multiplication by $\alpha_i(a)$. The $C^*$-correspondences $E_i \ot E_j$ may be
naturally identified with $_{\alpha_i \circ \alpha_j}\alg$; this provides immediately natural isomorphisms
$t_{i,j}$, as precisely described in \cite{Soltwo}. It is also explained there that  each representation of the
resulting product system $\be$ is uniquely determined by a  representation $\sigma$ of $\alg$ on some
Hilbert space $\Hil$ and a tuple of contractions $S_1,\ldots, S_k$ in $B(\Hil)$ satisfying the conditions
\begin{equation} \label{comm} \sigma(a) S_i = S_i \sigma(\alpha_i(a))\end{equation}
($i \in \onek, a \in \alg$); the representation is doubly commuting if and only if the tuple $S_1,\ldots, S_k$ is
doubly commuting. Similarly it may be shown that the corresponding representation is isometric if and only if each
$S_i$
is an isometry, and fully coisometric if and only if each $S_i$ is unitary.

The induced representation is determined by a  representation $\pi$ of $\alg$ on a Hilbert space
$\Kil$. The tuple  $S_1,\ldots, S_k$ forms a standard shift: $\Hil= l^2(\Nk) \ot \Kil$ and
\[ S_i (\xi \delta_m) = \xi \delta_{m+e_i}\]
$(i \in \onek, \xi \in \Kil, m\in \Nk$). What is nontrivial (and dependent on $\pi$) is the representation
$\sigma_{\pi}:\alg \to B(\Hil)$. It is given by the formula:
\[ \sigma_{\pi} (a)  (\xi \delta_m) = \pi (\alpha(m) (a)) \xi \delta_m,\]
where $a \in \alg, \xi \in \Kil, m\in \Nk$ and $\alpha(m) = \alpha_1^{m_1} \circ \cdots \circ \alpha_k^{m_k}$.

It follows from the above that the Wold decompositions from Theorems \ref{main} and \ref{main2} are equal (in
terms of respective Hilbert subspaces of the carrier space $\Hil$) to the standard Wold decompositions of the
tuple ($S_1,\ldots, S_k$) estalished in \cite{Sloc} and \cite{Popo}. The new aspect of Theorems  \ref{main} and
\ref{main2} reduces here to the identification of actions of $\alg$.

\section{Unitary extensions}

In the classical theory one of the first consequences of the Wold decomposition
is the fact that each isometry may be extended to a unitary (a one-sided shift
to a two-sided one). Independently it is also known that  every contraction can be extended to a coisometry.
In Section 5 of \cite{MSgen} P.\,Muhly and B.\,Solel obtain far reaching generalisations
of the latter result by exhibiting an explicit inductive extension procedure for a
completely contractive covariant representation of a $C^*$-correspondence
(see also \cite{MSproc}, where the particular case of a $C^*$-correspondence
induced by a unital injective endomorphism of $\alg$ is treated). Their construction uses the methods analogous
to these of the classical operator theory and under certain technical assumptions provides a fully coisometric extension
of a completely contractive representation. Apart from some special cases it is not known whether obtained in
such a way fully coisometric extension of an isometric representation remains isometric.

We are interested here in exploiting the Wold decomposition to obtain extensions which are simultaneously isometric
and fully coisometric, which we will further call \emph{unitary extensions}.
Our main result, Theorem \ref{unit1}, is very similar to Corollary 5.17 of \cite{MSgen}, only that the extension we
construct is simultaneously isometric and fully coisometric (i.e.\ unitary). Our approach is different to that of \cite{MSgen} and is based on the representation theory for $C^*$-algebras. The starting point is  Lemma 5.2,
 giving an easy representation theoretical criterion for the existence of coisometric extensions.
We then use Rieffel induction to establish the existence of a unitary extension under the assumption that the left action is faithful and has values in the algebra of compact operators. The section closes with the discussion of some examples
and a remark concerning the unresolved (and interesting) higher rank case.

To simplify the notation whenever $(\pi,\Kil)$ is a representation of $\alg$
the representation $(\phi \otimes \id_{\Kil}, E\ot_{\pi} \Kil)$ will be denoted
simply by $\phi \ot_{\pi} \id_{\Kil}$.

Let us fix the formal definition, stressing that we are assuming our extensions to be also isometric,
contrary to \cite{MSgen}:

\begin{deft}
Let $(\sigma,T)$ be an isometric representation of a $C^*$-correspondence
$E$ on a Hilbert space $\Hil$. We call an isometric representation $(\rho,V)$ of $E$ on a Hilbert space $\Kil$ containing $\Hil$ as a subspace
a unitary extension of $(\sigma,T)$ if the following conditions are satisfied:
\begin{rlist}
\item $(\rho,V)$ is fully coisometric;
\item for all $a \in \alg$, $\xi \in E$
and $h \in \Hil$
\[ V(\xi) h = T(\xi) h, \;\;\; \rho(a) h = \sigma(a) h.\]
\end{rlist}
\end{deft}

It is clear that when $\alg=\bc$ and $E=\bc$ the notion defined above is equivalent to the usual
concept of a unitary extension of an isometry (note that we require the extension to remain isometric).

Considering the question whether a given isometric representation $(\sigma, T)$ of $E$ has a unitary
extension, it is enough to focus on the `induced' part of the Wold decomposition. Suppose then that
$(\sigma, T)$ is induced from a  representation $(\pi, \Kil)$ of $\alg$. It is clear that the range of the isometric operator $\wT: E \ot \FFock(E) \ot_{\pi} \Kil \to \FFock(E) \ot_{\pi} \Kil$ is equal to
$E \ot \FFock(E) \ot_{\pi} \Kil =(\FFock(E) \ot_{\pi} \Kil) \ominus \Kil$. Suppose that $(\rho, V)$ is a unitary extension of $(\sigma,T)$. Then the representation space of $(\rho, V)$ decomposes as $(\FFock(E) \ot_{\pi} \Kil)
\oplus \Kil'$, and $\rho=  \sigma \oplus \pi'$. Moreover, we have the decomposition $\wV = \wT + \widetilde{V'}$, where
$\widetilde{V'}: E \ot \Kil' \to (\FFock(E) \ot_{\pi} \Kil) \oplus \Kil'$.
Since $(\rho, V)$ is fully coisometric, $\widetilde{V}$ is unitary and $\widetilde{V'}$ has to be an isometry
with the range equal to $\Kil \oplus \Kil'$. Using Lemma 2.1 of \cite{MSWold} a moment of thought allows to see that the sufficient and necessary
condition for $(\rho, V)$ to be a  unitary extension is that
$\widetilde{V'}: E \ot \Kil' \to \Kil \oplus \Kil'$ is a unitary intertwining the representations
$\phi \ot_{\pi'} \id_{\Kil'}$ and $\pi \oplus \pi'$.
 Let us summarise this discussion in the following lemma:

\begin{lem} \label{eqrep}
Let $(\pi, \Kil)$ be a  representation of $\alg$
and let $(\sigma,T)$ be an isometric representation of a $C^*$-correspondence
$E$ induced from $(\pi, \Kil)$.
The representation $(\sigma,T)$ has a unitary extension if and only if there is
a  representation $(\pi', \Kil')$ of $\alg$ such that
\begin{equation} \label{plus} \phi \ot_{\pi'} \id_{\Kil'} \approx \pi \oplus \pi'.\end{equation}
\end{lem}

Note that if $\alg=\bc$, the equivalence of different representations of $\alg$ amounts to checking that
the Hilbertian dimensions of relevant Hilbert spaces are equal. In general the problem of constructing, for a given $\pi$, a representation $\pi'$ such that \eqref{plus} holds is clearly more complicated.

To obtain some general statements it is necessary to assume that the left action of $\alg$ on $E$, denoted further by $\phi$, is faithful.
Indeed, consider the following
example. Let $\alg=\bc^2:= \bc e_1 \oplus \bc e_2$, $E=\alg$ as a $C^*$-Hilbert module (so that the scalar product is given via the formula $\la a, b \ra = a^*b$ and the right action of $\alg$ by the right multiplication), with the
left action given by $\phi(\lambda e_1 + \mu e_2) (a)= \lambda a$ ($a\in \alg, \lambda, \mu \in \bc$). Note that $\phi$ is not faithful.
If $(\pi,\Kil)$ is a  representation of $\alg$, it decomposes as $(\pi_1 \oplus \pi_2, \Kil_1 \oplus \Kil_2)$,
where $\pi_1 (e_1) = I_{\Kil_1}$, $\pi_2 (e_2) = I_{\Kil_2}$, $\pi_1(e_2) = \pi_2(e_1)=0$.
Consider the representations $\sigma_1:=\phi \ot_{\pi_1} \id_{\Kil_1}$,
$\sigma_2:=\phi \ot_{\pi_2} \id_{\Kil_2}$. It is easy to check that $E \ot_{\pi_1} \Kil_1 \approx \bc \ot \Kil_1$,
$E \ot_{\pi_2} \Kil_2 \approx \bc \ot \Kil_2$ as Hilbert spaces. Further
$\phi \ot_{\sigma_1} \id_{\Kil_1} \approx \pi_1$, but $\phi \ot_{\sigma_2} \id_{\Kil_2} (e_2)=0$, $\phi \ot_{\sigma_2} \id_{\Kil_2} (e_1) = I_{\Kil_2}$.
This implies that the representation of $\alg$ arising as $\phi \ot \id$ is always zero on $e_2$. If then $\pi$ is given as above and $\Kil_2$ is nonzero,
then there cannot exist any representation $\pi'$ of $\alg$ such
that \eqref{plus} is satisfied.

The $C^*$-algebra of compact operators on a Hilbert module $E$ will be denoted by $\comp$
(see \cite{Lance}).
The following lemma is based on the standard techniques from  Morita equivalence theory (see \cite{Morita}). In order to avoid
the representations to become too large we need some mild bounds of the size of the constructions. For a normed linear space $X$ denote
by $w(X)$ the smallest cardinal of a total subset. If $\Hil$ is a Hilbert space then $w(\Hil)$ equals the Hilbertian dimension of $\Hil$.
It is easy to see that $w(E \ot_{\rho} \Hil) \leq w(E)w(\Hil)$ for any representation $(\rho,\Hil)$.  Notice also that $w(\comp) \leq w(E)^2$.

\begin{lem}        \label{Morita}
Suppose that $E$ is a $C^*$-correspondence with faithful left action of $\alg$ on  $E$
and the range of $\phi$ contained in the space of compact operators of $E$,
in other words, $\phi: \alg \to \comp$ is an injective $*$-homomorphism.
Then given any representation $(\rho ,\Hil)$ of $A$ there is a representation
$(\hat{\rho},\hat{\Hil})$ of $A$ such that $(\phi \otimes_{\hat{\rho}} \id_{\hat{\Hil}},\Kil)$ contains a subrepresentation
equivalent to $\rho$. Moreover, one can arrange that $w(\hat{\Hil}) \leq w(E)^3 w(\Hil)$ and $w(\Kil) \leq w(E)^2 w(\Hil)$.
\end{lem}

\begin{proof}
Let $\rho$ be non-zero.
The map $\rho \circ
\phi^{-1}$ is a representation of the $C^*$-algebra $\phi(\alg)$.
As $\phi(\alg)$ is a $C^*$-subalgebra of $\comp$,  $\rho \circ \phi^{-1}$ can be extended
(Proposition II.6.4.11, \cite{Black})
to a  representation  $\sigma$ of $\comp$, possibly on
a bigger Hilbert space $\Kil \supseteq \Hil$.  The precise statement is that there is a
representation $\sigma$ of $\comp$ such that $\sigma|_{\phi(\alg)}$ contains
 $\rho \circ \phi^{-1}$ as a subrepresentation. The standard construction uses extensions of states and
the GNS construction, which shows that one can assume that $w(\Kil) \leq w(\comp) w(\Hil) \leq w(E)^2 w(H)$.
From the theory of Rieffel correspondence it follows that $\sigma$ must arise from
some representation $(\hat{\rho}, \hat{\Hil})$ of $\alg$ via tensoring with the identity on $E$.
Let us briefly recall the argument for the reader's convenience: the dual module $\wt{E}$
is an $I-\comp$ equivalence module, where $I$ is the ideal generated by all inner products in $E$.
The induction $(\sigma^{\wt{E}}, \wt{E} \ot_{\comp} \Kil )$ and the left action $\wt{\phi}(a) \wt{\xi} =  \wt{\xi a^*}$
provide a representation $\hat{\rho} : A \to B (\hat{\Hil} )$, where $\hat{\Hil}= \wt{E} \ot_{\comp} \Kil$.
Using $E \ot_A \wt{E} = E \ot_{I} \wt{E} = \comp$ it follows that
the representation ${\hat{\rho}}^E : \adj \to B(E \ot_{\hat{\rho}} \hat{\Hil})= B(\Kil)$ is the
canonical extension of $\sigma$ (up to unitary equivalence).
Thus $\phi\otimes_{\hat{\rho}} \id_{\hat{\Hil}}$ is
equivalent to $\sigma \circ \phi$. By our construction $\sigma \circ \phi$ contains $\rho$ as a
subrepresentation. Finally  $w(\hat{\Hil}) \leq w(\wt{E})w(E)^2 w(\Hil)= w(E)^3 w(\Hil)$ is clear from the construction.
\end{proof}

We do not know whether the assumption $\phi(\alg) \subset \comp$ may be weakened.
The problem  lies in the fact that Morita equivalence provides a correspondence between
representations of $\alg$ and $\comp$. In general, when we only know that $\phi$ has values in $\adj$
(the algebra of adjointable operators), the representation $\rho \circ \phi^{-1}$ from the proof above
may still be extended to a representation of $\adj$, but it may well happen that the
restriction of this extension to $\comp$ is trivial (see the example after Theorem \ref{unit1}).

Recall that the universal representation $(\pi_u, \Hil_u)$ of $\alg$ is the direct sum over `all' representations
of $\alg$. Avoiding set theoretic difficulties we define it as the direct sum over all representations of $\alg$
on a fixed Hilbert space $\Hil_0$ with $w(\Hil_0)= w(\alg)$.

\begin{tw} \label{unit1}
Suppose that $E$ is a $C^*$-correspondence with faithful left action of $\alg$ on  $E$
and range contained in the space of compact operators of $E$.
Then every isometric representation of
$E$ on a Hilbert space $\Hil$ has a unitary extension.
\end{tw}
\begin{proof}
Consider an isometric representation $(\sigma, T)$ of $E$. We may and do assume that
$(\sigma, T)$ is induced from a representation $(\pi, \Kil)$ of $\alg$.
From Lemma \ref{eqrep} it follows that it suffices to find a
representation $(\pi', \Kil')$ of $\alg$ such that
$\phi \ot_{\pi'} \id_{\Kil'} \approx \pi \oplus \pi'$.
By Lemma \ref{Morita} applied to the the universal representation $\pi_u$ we find a representation $(\hat{\pi}_u, \hat{\Hil}_u)$
such that $\pi_u \preceq \phi \ot_{\hat{\pi}_u}  \id_{\hat{\Hil}_u}$ and every (non-zero) subrepresentation
of $\hat{\pi}_u$ has multiplicity not exceeding $w(\Hil_u)w(E)^3$. For any cardinal $\aleph$
denote by $\pi_u^{\aleph}$ the direct sum of $\aleph$ copies of $\pi_u$. Put $\aleph= w(\Hil_u)w(E)^3$. Then it follows that
$\phi \ot_{\hat{\pi}_u} \ot \id_{\Hil_u} \preceq \phi \ot_{\pi_u^{\aleph}} \id_{\Hil_u^{\aleph}}$.
Moreover, by the remark before Lemma \ref{Morita}, $\phi \ot_{\pi_u^{\aleph}} \id_{\Hil_u^{\aleph}} \preceq  \pi_u^{w(E) w(\Hil_u) \aleph}$.
Thus letting $\aleph' = w(E) w(\Hil_u) \aleph_0$ we have
$$
\pi_u^{\aleph'} \preceq \phi \ot_{\pi_u^{\aleph'}}  \id_{\Hil_u^{\aleph'}} \preceq \pi_u^{\aleph'} \approx \pi \oplus \pi_u^{\aleph'}.
$$
Since the multiplicity of any subrepresentation of $\pi_u^{\aleph'}$ and $\phi \ot_{\pi_u^{\aleph'}} \id$ is $\aleph'$ it follows that they are unitarily equivalent
so that $\pi' = \pi_u^{\aleph'}$ is as required.
\end{proof}
Note that if $\pi$ is a type I representation, one can avoid most of the cardinality considerations. It is enough then to use the fact that $\pi_a \leq \phi \ot_{\pi_a} \id_{\Hil_a}$,
where $(\pi_a, \Hil_a)$ denotes the reduced atomic representation of $\alg$. The latter fact follows from the observation that Rieffel correspondence behaves well with respect to the irreducibility of representations involved.

As already pointed out the result above is closely related to Corollary 5.17 of \cite{MSgen}.
The methods used in \cite{MSgen} are different
to these used here: the fully coisometric extension is obtained via an inductive procedure, in a way analogous to that of constructiong isometric dilations. It is possible that when $(T, \sigma)$ is isometric, the construction given by P.\,Muhly and B.\,Solel yields in fact an isometric representation. To show it, it would be sufficient to establish that the isometric property is preserved
under each step of the inductive reasoning (i.e.\ the step given in Proposition 5.7 of \cite{MSgen}). This is the case in the simplified framework of $C^*$-correspondences considered in \cite{MSproc}. In general it remains an open question.

As Example 5.16 of \cite{MSgen} shows, it may happen that a completely contractive covariant representation $(T, \sigma)$ of $E$ has no fully
coisometric extension. Although the representation  considered there is not isometric, the example may be easily modified so that we obtain an isometric representation $(\sigma,T)$ of a $C^*$-correspondence $E$, with $E$ nondegenerate and the left action $\phi:\alg \to \adj$ faithful, which has no fully
coisometric extension. We recall it below for completeness:

\begin{exam}[\cite{MSgen}]
Let $\alg = l^{\infty}$, and let $E=C^{\infty}(\alg)$ denote the column module over $\alg$ (the space of all sequences $(a_i)_{i \in \bn}$ of elements of $\alg$ such that $\sum_{i=1}^{\infty} a_i ^* a_i$ converges in $\alg$) equipped with the natural $\alg$-valued scalar product and the left action given by entry-wise multiplication. It may be checked that $E$ is nondegenerate, the left action $\phi$ is faithful,
$\phi(c_0) \subset \comp$ and $\phi(c_o) E$ is dense in $E$.
Suppose that $(\pi, \Kil)$ is a nonzero nondegenerate representation of $\alg$ such that
$\pi|_{c_0} =0$ and let $(\sigma, T)$ be an isometric representation of $E$ induced by $(\pi, \Kil)$.
Then as the intersection of $\sigma(c_0) (\FFock(\be) \ot_{\pi} \Kil)$ with $\Kil$ is trivial, Proposition 5.13 of \cite{MSgen} implies that $(\sigma, T)$ cannot have a fully coisometric extension.
\end{exam}

Theorem \ref{unit1}, although quite general, does not cover the following particular case:
it is clear from the discussion before Lemma \ref{eqrep} that if $\alg=\bc$ then every
isometric representation of $E$ has a unitary extension and at the same time when $E$ is
not finite dimensional one cannot expect that $\phi$ acts on $E$ by compact operators. To cover
 this case we establish the following:

 \begin{tw}   \label{unit2}
Let $(\sigma,T)$ be an isometric representation of a $C^*$-correspondence
$E$ on a Hilbert space $\Hil$. Suppose that the induced part of $(\sigma,T)$
arises from a representation $(\pi, \Kil)$ of $\alg$.
If $(\pi, \Kil)$ is equivalent to a subrepresentation of
$\phi \ot_{\pi} \id_{\Kil}$
then $(\sigma,T)$ has a unitary extension.
\end{tw}

\begin{proof}
Again, we may and do assume that
$(\sigma, T)$ has trivial fully coisometric part. Consider the representation $\sigma^{\infty}$
given by the countably infinite number of copies of $\sigma$. As by our assumption
\[ \pi^{\infty} \oplus (\phi \ot_{\pi} \id_{\Kil})^{\infty} \approx (\phi \ot_{\pi} \id_{\Kil})^{\infty},\]
there is also
\[ \phi \ot_{\sigma^{\infty}} \id_{\Hil^{\infty}} \approx (\phi \ot_{\sigma} \id_{\Hil})^{\infty}
\approx \pi^{\infty} \oplus (\phi \ot_{\pi} \id_{\Kil})^{\infty} \oplus \cdots \approx
(\phi \ot_{\pi} \id_{\Kil})^{\infty} \oplus \cdots \approx \sigma^{\infty}.\]
Finally
\[\pi \oplus \sigma^{\infty} \approx \pi \oplus \pi^{\infty} \oplus (\phi \ot_{\pi} \id_{\Kil})^{\infty} \oplus \cdots
\approx \pi^{\infty} \oplus (\phi \ot_{\pi} \id_{\Kil})^{\infty} \oplus \cdots = \sigma^{\infty},\]
and Lemma \ref{eqrep} ends the proof.
\end{proof}

Assumptions of the above theorem are trivially satisfied for each $(\sigma, T)$
if $\alg= \bc$, as then the action of $\phi$ is given by an orthogonal projection in $E$, and
$\phi \ot_{\pi} \id_{\Kil}$ is equivalent to a certain number of copies of $\pi$.
In particular for $\alg= \bc$, $E = \bc^n$ ($n \in \bn$) we obtain another proof of
Proposition 2.6 in \cite{Popescu}. Note that in fact in the above mentioned paper there is
potential for a certain confusion in terminology.
If the notion of dilation is understood to include certain co-invariance conditions (as it is for
example in this paper, see the definition before Proposition \ref{dilatafter}), it is in general
impossible to dilate a contraction to a unitary: we first \emph{dilate} a contraction to an isometry
and then \emph{extend} the resulting isometry to a unitary. We leave the formal deconstruction of the
remark above to the reader.

A modification of the example presented after Theorem \ref{unit1} shows
that it may happen that the assumptions of Theorem \ref{unit1} are satisfied, whereas those of Theorem \ref{unit2}
are not. Consider again $\alg = \bc^2$, $E= \bc^2$ and let $\phi$  be now given
by the formula
\[ \phi(\lambda e_1 + \mu e_2) \binom{\lambda'}{\mu'} = \binom{\mu\lambda'}{\lambda\mu'}, \;\;\; \lambda, \lambda', \mu, \mu' \in \bc.\]
Let $(\pi, \Kil)$ be a representation of $\bc^2$ given by
$\pi(e_1) = I_{\Kil}$, $\pi (e_2) = 0$. Following the same arguments as before one can see that $\phi \otimes_{\pi} \id_{\Kil}$ is  equivalent to  $(\pi', \Kil)$, where
$\pi'(e_1) = 0$, $\pi' (e_2) = I_{\Kil}$. The representations $\pi$ and $\phi \otimes_{\pi} \id_{\Kil}$
are therefore disjoint.
On the other hand it is easy to check that if $(\pi_u, \Hil_u)$ is the universal representation
of $\bc^2$ then  $\phi \otimes_{\pi_u} \id_{\Hil_u} \approx \pi_u$. As $E$ is finitely generated,
$\adj=\comp$ and the range of $\phi$ is obviously contained in $\comp$.

\subsection*{Minimality and difficulties in the case $k>1$}
 In the classical context of a
unitary extension $U\in B(\Kil)$ of the isometry $T\in B(\Hil)$, the minimal closed subspace of $\Kil$ containing $\Hil$ and
invariant under both $U$ and $U^*$ is generated by $\{U^m h: m \in \bz, h \in \Hil\}$. The restriction of $U$ to his subspace
provides the minimal unitary extension which is unique up to unitary equivalence.

In the framework described above there is a similar concept of minimality. Given a unitary extension $(\pi,V)$ on $\Kil$
of an isometric representation of $E$ on $\Hil$ one can find minimal reducing subspaces for $(\pi,V)$ such that the corresponding restriction is still a unitary extension of the original representation. Indeed, $(\pi,V)$ corresponds to a representation of
the Cuntz-Pimsner algebra $\Oo_E$ and the minimal extension is given by the cyclic subspace of this representation generated by $\Hil$. Uniqueness (unless e.g. $\alg = \bc$) is however in general lost (we refer to \cite{MSproc} for the discussion of the reasons behind it). There is also no reason for the constructions given in Theorems \ref{unit1} and \ref{unit2} to yield minimal extensions.


In the classical context the minimal unitary extension provides a starting point for the It\^o construction of a
joint extension of a tuple of commuting isometries to a tuple of commuting unitaries (Proposition 6.2, \cite{Nagy}).
If however $\textrm{card} \mathcal{J} >1$, the crucial inductive step allowing to establish
that a natural extension of one row-isometry to the minimal `extension space' of another commuting row-isometry
remains a row-isometry, fails. It would be therefore highly nontrivial and very interesting
to establish some counterparts of  Theorems
\ref{unit1} and \ref{unit2} for the case $k>1$.

\subsection*{ACKNOWLEDGMENT}
We would like to thank the anonymous referee for the careful reading of our manuscript and for bringing into our attention the connections between our work on coisometric extensions and
the results of Section 5 of \cite{MSgen} and of \cite{MSproc}.


\begin{thebibliography}{MSWold}

\bibitem [BCL] {Coburn} C.A.\,Berger, L.A.\,Coburn and A.\,Lebow,
Representation and index theory for $C^*$-algebras generated by commuting isometries,
  \emph{J. Funct. Anal.}  \textbf{27} (1978) 51--99.

\bibitem [Bla] {Black} B.\,Blackadar, ``Operator Algebras. Theory of $C^*$-Algebras and
Von Neumann Algebras", Encyclopaedia of Mathematical Sciences, 122. Operator Algebras and Non-commutative Geometry, III. Springer-Verlag, Berlin, 2006.

\bibitem [Fow] {fowl} N.\,Fowler, Discrete product systems of Hilbert bimodules,
  \emph{Pacific J. Math.}  \textbf{204} (2002) 335--375.

\bibitem [JuK] {jury}  M.T.\,Jury and D.W.\,Kribs,  Ideal structure in free semigroupoid algebras from directed graphs, \emph{J. Operator Theory} \textbf{53}  (2005),  no. 2, 273--302.

\bibitem [Lan] {Lance} E.C.\, Lance, ``Hilbert $C^*$-modules", LMS Lecture Notes Series \textbf{210},
Cambridge University Press, Cambridge, 1995.

\bibitem [MS$_1$] {MSgen} P.\,Muhly and B.\,Solel, Tensor algebras over $C^*$-correspondences (Representations, dilations and $C^*$-envelopes),  \emph{J. Funct. Anal.}  \textbf{158} (1998) 389--457.

\bibitem [MS$_2$] {MSWold} P.\,Muhly and B.\,Solel, Tensor algebras, induced representations, and the Wold decomposition,   \emph{Canad. J. Math.}  \textbf{51}  (1999),  no. 4, 850--880.

\bibitem [MS$_3$] {MSproc} P.\,Muhly and B.\,Solel, Extensions and dilations for $C^*$-dynamical systems,  \emph{Operator theory, operator algebras, and applications},  375--381, \emph{Contemp. Math.}, \textbf{414}, Amer. Math. Soc., Providence, RI, 2006.

\bibitem [Nic] {Nica} A.\,Nica, $C^*$-algebras generated by isometries and Wiener-Hopf operators, \emph{J.\,Operator Theory} \textbf{27}  (1992),  no. 1, 17--52.

\bibitem [Pope] {Popescu} G.\,Popescu, Isometric dilations
for infinite sequences of noncommuting operators,
  \emph{Trans.\,AMS}  \textbf{316} (1989), no.\,2, 523--536.

\bibitem [Pop] {Popo} D.\,Popovici, A Wold-type decomposition for commuting isometric pairs,
  \emph{Proc.\,AMS}  \textbf{132} (2004), no.\,8, 2303--2314.

\bibitem [RaW] {Morita} I.\,Raeburn and D.\,Williams, ``Morita Equivalence and Continuous-Trace-$C^*$-Algebras", Mathematical Surveys and Monographs, 60. American Mathematical Society, Providence, RI, 1998.

\bibitem[S\l o] {Sloc} M.\,S\l oci\'nski, On the Wold-type decomposition of a pair of commuting isometries,
\emph{Ann. Polon. Math.} \textbf{37}  (1980), no.\,3, 255--262.

\bibitem [So$_1$] {Soltwo} B.\,Solel, Representations of product systems over semigroups and dilations of commuting CP maps,  \emph{J. Funct. Anal.}  \textbf{235}  (2006),  no. 2, 593--618.

\bibitem [So$_2$] {Solk} B.\,Solel, Regular dilations of representations of product systems, \emph{preprint}, math.OA/0504129.

\bibitem [SzF] {Nagy} B.\,Sz.-Nagy and C.\,Foias, ``Harmonic analysis of operators on Hilbert space", North Holland, Amsterdam, 1970.

\end{thebibliography}
\end{document}